\documentclass[english, 12pt,leqno]{amsart}

\usepackage{amsmath}
\usepackage{amsfonts}
\usepackage{amssymb}

\newtheorem{Def}{D\'efinition}

\newtheorem{Cor}[Def]{Corollary}
\newtheorem{Exam}[Def]{Example}
\newtheorem{Lemma}[Def]{Lemme}
\newtheorem{Prop}[Def]{Proposition}

\newtheorem{Rem}[Def]{Remark}
\newtheorem{Thm}[Def]{Theorem}

\title[The Fuglede-Kadison determinant]
{The Fuglede-Kadison determinant
\\
theme and variations}

\author[Pierre de la Harpe]
{Pierre de la Harpe}

\date{July 6, 2011. Minor corrections on July 7, 2012}

\address{Pierre de la Harpe,
Section de math\'ematiques, 
Universit\'e de Gen\`eve, 
C.P.~64, 
CH--1211 Gen\`eve 4. 
\emph{E-mail address:} Pierre.delaHarpe@unige.ch
}



\begin{document}

\begin{abstract}
We review the definition of determinants for finite von Neumann algebras,
due to Fuglede and Kadison (1952), 
and a generalisation 
for appropriate groups of invertible elements in Banach algebras, 
from a paper by Skandalis and the author (1984).
After some reminder on K-theory and Whitehead torsion,
we hint at the relevance of these determinants
to the study of $L^2$-torsion in topology. Contents:
\begin{itemize}
\item[\ref{classical}.]
The classical setting.
\item[\ref{vna}.]
On von Neumann algebras and traces.
\item[\ref{FK1}.]
The Fuglede-Kadison determinant for finite von Neumann algebras.
\item[\ref{motivating}.]
A motivating question.
\item[\ref{K}.]
A reminder on $K$-theory.
\item[\ref{FK2}.]
Revisiting the Fuglede-Kadison and other determinants.
\item[\ref{baratintorsion}.]
On Whitehead torsion.
\item[\ref{baratinL2torsion}.]
A few lines  on $L^2$-torsion.
\end{itemize}
\end{abstract}

\maketitle

\section{\textbf{The classical setting}}
\label{classical}

\subsection{Determinants of matrices over commutative rings}
Let $\mathcal R$ be a  ring with unit.
For an integer $n \ge 1$, denote by $\operatorname{M}_n(\mathcal R)$
the ring of $n$-by-$n$ matrices over $\mathcal R$
and by $\operatorname{GL}_n(\mathcal R)$  its group of units. 
$\mathcal R^*$ stands for $\operatorname{GL}_1(\mathcal R)$.
\par

\emph{Suppose $\mathcal R$ is commutative.}
The determinant 
\begin{equation}
\label{eqdetalgclassique}
\det \, : \, \operatorname{M}_n(\mathcal R) \longrightarrow \mathcal R
\end{equation}
is defined by a well-known explicit formula, polynomial in the matrix entries.
It is alternate multilinear in the columns of the matrix,
and normalised by $\det (1_n) = 1$;
when $\mathcal R$ is a field, these properties constitute
an equivalent  definition, 
as was lectured on by Weierstrass and Kronecker probably in the 1860's,
and published much later 
(see \cite[Siebzehnte Vorlesung, Page 291;
published 12 years after Kronecker's death]{Kron--03}).
\par

For  $x,y \in \operatorname{M}_n(\mathcal R)$,
we have $\det (xy) = \det (x) \det (y)$.
For ${x \in  \operatorname{M}_n(\mathcal R)}$ with $\det (x)$ invertible,
an explicit formula shows that $x$ itself is invertible,
so that $\det (x) \in \mathcal R ^*$ 
if and only if $x \in \operatorname{GL}_n(\mathcal R)$.
The restriction
\begin{equation}
\label{eqdetgrpeclassique}
\operatorname{GL}_n(\mathcal R) \longrightarrow \mathcal R^*, 
\hskip.5cm x \longmapsto \det x
\end{equation}
is a group homomorphism.
\medskip

\subsection{Three formulas for complex matrices involving
determinants, exponentials, traces, and logarithms}
\label{threeformulas}
Suppose that $\mathcal R$ is \emph{the field $\mathbf C$ of complex numbers.}
The basic property of determinants that we wish to point out
is the relation
\begin{equation}
\label{detexpC}
\det (\exp y) \, = \, \exp( \operatorname{trace} (y))
\hskip.5cm \text{for all} \hskip.2cm y \in \operatorname{M}_n(\mathbf C) .
\end{equation}
Some exposition books give this as a very basic formula \cite[$\S$~16]{Arno--73};
it will reappear below as Equation (\ref{detexpA}).
It can also be written
\begin{equation}
\label{detx}
\det (x) \, = \,  \exp( \operatorname{trace} (\log x))
\hskip.5cm \text{for appropriate} \hskip.2cm 
x \in \operatorname{GL}_n(\mathbf C) .
\end{equation}
``Appropriate''  can mean several things.
If  $\Vert x - 1 \Vert < 1$,
then $\log x$ can be defined by the convergent series
\begin{equation*}
\log x \, = \, \log \big( 1 + (x-1) \big) \, ) \, = \, 
\sum_{k=1}^\infty \frac{(-1)^{k-1}}{k} (x-1)^k .
\end{equation*}
If $x$ is conjugate to a diagonal matrix,
then $\log x$ can be defined compo\-nent-wise
(in pedantic terms, this is functional calculus,
justified by the spectral theorem).
In (\ref{detx}), note that
the indeterminacy in the choice of the logarithm of a complex number
is swallowed by the exponential, because $\exp 2 \pi i = 1$.
\par

Let $x \in \operatorname{GL}_n(\mathbf C)$.
Since the group is connected, we can choose a piecewise smooth path
$\xi : [0,1] \longrightarrow \operatorname{GL}_n(\mathbf C)$ 
from $1$ to $x$.
Since $\log \xi (\alpha)$ is a primitive of $\dot{\xi} (\alpha) \xi (\alpha)^{-1} d\alpha$,
it follows from (\ref{detx}) that
\begin{equation}
\label{detexpC2}
\det (x) \, \overset{!}{=} \, 
\exp \int_0^1  \operatorname{trace}(\dot{\xi} (\alpha) \xi (\alpha) ^{-1}) d\alpha .
\end{equation}
This will be our motivating formula for Section~\ref{FK2}, 
and in particular for Equation (\ref{Eqbasic}).
\par

The sign $\overset{!}{=}$ stands for a genuine equality,
but indicates that some comment is in order.
A priori, the integral depends on the choice of $\xi$,
and we have also to worry about the determination of $\log \xi (\alpha)$.
As there is \emph{locally} no obstruction to choose a continuous
determination of the primitive 
$\log \xi (\alpha)$ of $\dot{\xi} (\alpha) \xi (\alpha)^{-1} d\alpha$,
the integral is invariant by small changes of the path (with fixed ends),
and therefore depends only on the homotopy class of $\xi$,
so that it is defined modulo its values
on (homotopy classes of) closed loops.
The fundamental group $\pi_1(\operatorname{GL}_n(\mathbf C))$
is infinite cyclic, generated by the homotopy class of
\begin{equation*}
\xi_0 \, : \,  [0,1] \longrightarrow \operatorname{GL}_n(\mathbf C),
\hskip.5cm
\alpha \longmapsto
\begin{pmatrix}
e^{2 \pi i \alpha} & 0 \\ 0 & 1_{n-1}
\end{pmatrix} ,
\end{equation*}
and we have
$\int_0^1  \operatorname{trace}(\dot{\xi_0} (\alpha) \xi_0 (\alpha) ^{-1}) d\alpha
= 2 \pi i$.
Consequently, the integral in the right-hand side of (\ref{detexpC2})
is defined modulo $2 \pi i \mathbf Z$, 
so that the right-hand side itself is well-defined.
(This will be repeated in the proof of Lemma \ref{lemmepourtildedelta}.)
\par

Since a connected group is generated by any neighbourhood of the identity,
there exist $x_1, \hdots, x_k \in \operatorname{GL}_n(\mathbf C)$
such that $x = x_1 \cdots x_k$
and $\Vert x_j - 1 \Vert < 1$ for $j = 1, \hdots, k$,
and  one can choose
\begin{equation*}
\xi(\alpha) \,  = \,  \exp( \alpha (\log x_1) ) \cdots \exp( \alpha (\log x_k) ) .
\end{equation*}
A short computation with this $\xi$ gives
\begin{equation*}
\exp \int_0^1  \operatorname{trace}(\dot{\xi} (\alpha) \xi (\alpha) ^{-1}) d\alpha
\, = \, \exp ( \operatorname{trace}( \log x_1 )) \cdots
\exp ( \operatorname{trace}( \log x_k ))
\end{equation*}
and it is now obvious that (\ref{detx}) implies (\ref{detexpC2}).

\subsection{Historical note}
Determinants arise naturally with \emph{linear systems of equations},
first with $\mathcal R = \mathbf R$, 
and more recently also with $\mathcal R = \mathbf C$.
They have a prehistory
in Chinese mathematics from the 2nd century BC \cite{MacT}.
In modern Europa, there has been an early contribution by\footnote{There are
also resultants and determinants in the work of the
Japanese mathematician Seki Takakazu,
a contemporary of Leibniz and Newton.}
Leibniz in 1693, unpublished until 1850.
Gabriel Cramer wrote an influential book, published in 1750.
Major mathematicians who have written about determinants include 
B\'ezout, Vandermonde, Laplace, Lagrange, 
Gauss, Cauchy, Jacobi, Sylvester, Cayley, ...
The connection between determinants of matrices in $M_3(\mathbf R)$
and volumes of parallelepipeds is often attribued to Lagrange (1773).
Let us  mention 
an amazing  book on the history of determinants 
\cite{Muir--23}:
four volumes, altogether more than 2000 pages,
an ancestor of  the
\emph{Mathematical Reviews}, for \emph{one} subject,
covering the period 1693--1900.
\par

There is an extension of (\ref{eqdetalgclassique}) to a skew-field $k$
by Dieudonn\'e, 
where the range of the mapping 
$\operatorname{det}$ defined on $\operatorname{M}_n(k)$
is\footnote{We denote by $D\Gamma$  the \emph{group of commutators}
of a group $\Gamma$.} 
$(k^*/Dk^*) \sqcup \{0\}$
(see  \cite{Dieu--43}, \cite{Arti--57},
and \cite{Asla--96} for a discussion when $k$
is the skew-field of Hamilton quaternions).
Determinants in case of a non-commutative ring $\mathcal R$
has motivated a lot of work,
in particular by Gelfand and co-autors since the early 1990's
\cite{GGRV--05}.
Let us also mention a version for super-mathematics due to Berezin
(see  \cite{Bere--79} and \cite[Chapter 3]{Mani--88}),
as well as ``quantum determinants'', of interest in low-dimensional topology
(see for example \cite{HuLe--05}).
\par

The notion of determinants extends to matrices over a ring without unit
(by ``adjoining a unit to the ring'').
In particular, in functional analysis, there is a standard notion of determinants
which appears in the theory of Fredholm integral equations,
for example for operators on a Hilbert space of the form $1+x$
where $x$ is ``trace-class'' \cite{Grot--56, Simo--79}.
\par

The oldest occurence I know of  $\exp y$ or $\log x$, \emph{including the notation},
defined by the familiar power series in the matrix $y$ or $x-1$,
is in \cite[Page 374]{Metz--92}; see also \cite{vNeu--29}.
But exponentials of linear differential operators appear also early in Lie theory,
see for example \cite[Page 75]{LieE--88} and  \cite[Page 82]{Hawk--00}, 
even if Lie never uses a notation like $\exp X$
(unlike Poincar\'e, see his $e^{\alpha X}$ in \cite[Page 177]{Poin--99}).
\par

There is a related and rather old formula known as the
 ``Abel-Liouville-Jacobi-Ostrogradskii Identity''.
Consider a homogeneous linear differential equation of the first order
$y'(t) = A(t)y(t)$, for an unknown function $y : [t_0,t_1] \longrightarrow \mathbf R^n$.
The columns of a set of $n$ linearly independent solutions constitute
the \emph{Wronskian matrix} $W(t)$.
It is quite elementary (at least nowadays!) to show that
$W'(t) = A(t)W(t)$, 
hence $\frac{d}{dt} \det W(t) = \operatorname{trace}(A(t)) \det W(t)$,
and therefore 
\begin{equation*}
\det W(t) = \det W(t_0) \exp\Big( \int_{t_0}^t \operatorname{trace}(A(s)) ds \Big) ,
\end{equation*}
a close cousin of Equation (\ref{detexpC2}).
The name of this identity refers to Abel (1827, case $n=2$), 
Liouville \cite{Liou--38},
Ostrogradskii (also 1838), and Jacobi (1845).
This was pointed out to me by Gerhard Wanner \cite[Section I.11]{HaNW--93};
also, Philippe Henry showed me this identity 
on the last but five line of \cite{Darb--80}
(which does not contain references to previous work).

\medskip

Finally, a few words about the authors of the 1952 paper alluded to in our title.
\emph{Bent Fuglede} is a Danish mathematician born in 1925.
He has been working on mathematical analysis;
he is also known for a book on
\emph{Harmonic maps between Riemannian polyhedra}
(co-author Jim Eells, preface by Misha Gromov).
\emph{Richard Kadison} is an
American mathematician, born in this same year 1925.
He is known for his many contributions to operator algebras;
his ``global vision of the field was certainly essential for my own development''
(words of Alain Connes,
when Kadison was awarded the Steele Prize 1999 for Lifetime Achievement,
see \cite{Kadi}).

\subsection{Plan}
Section \ref{vna} is a reminder on von Neumann algebras
based on three types of examples, 
Section \ref{FK1} is an exposition of the original Fuglede-Kadison idea,
Section \ref{motivating} stresses the difference between 
the complex-valued standard determinant
and the real-valued Fuglede-Kadison determinant,
Section \ref{K} is a reminder on some notions of K-theory.
Section \ref{FK2} exposes the main variations of our title: 
determinants defined for connected groups of invertible elements 
in complex Banach algebras.
We end by recalling in Section \ref{baratintorsion} 
a few facts about Whitehead torsion,
with values in $\operatorname{Wh}(\Gamma)$, 
which is a quotient of the group $K_1$ of a group algebra
$\mathbf Z [\Gamma]$,
and by alluding in Section \ref{baratinL2torsion}
to $L^2$-torsion, which is defined in terms of (a variant of) the Fuglede-Kadison determinant.

\subsection{Thanks}
I am grateful to 
Georges Skandalis for \cite{HaS--84a},
to Tatiana Nagnibeda and Stanislas Smirnov for their invitation
to deliver a talk in Saint Petersburg on this subject,
to Dick Kadison for encouragement to clean up my notes, 
as well as to Thierry Giordano, Jean-Claude Hausmann, 
Wolfgang L\"uck, Thierry Vust, and Claude Weber
for useful comments.

\section{\textbf{On von Neumann algebras and traces}}
\label{vna}


In a series of papers from 1936 to 1949, 
Francis Joseph Murray and John von Neumann
founded the theory of \emph{von Neumann algebras}
(in their terminology \emph{rings of operators}),
which are complex $*$-algebras
representable by unital weakly-closed $*$-subalgebras 
of some $\mathcal L (\mathcal H)$,
the algebra of all bounded operators on a complex Hilbert space $\mathcal H$.
\par

We will first give three examples of pairs $(\mathcal N, \tau)$,
with $\mathcal N$ a finite von Neumann algebra and $\tau$
a finite trace on it. We will then recall some general facts,
and define a few terms, such as ``finite von Neumann algebra'', ``finite trace'', 
and  ``factor of type $\operatorname{II}_1$''.

\begin{Exam}[\textbf{factors of type $\operatorname{I}_n$}]
\label{I_nvna}
For any $n \ge 1$, the matrix algebra $\operatorname{M}_n(\mathbf C)$
is a finite von Neumann algebra known as a \emph{factor of type $\operatorname{I}_n$}.
The involution is given by $(x^*)_{j,k} = \overline{x_{k,j}}$.
The linear form
$x \longmapsto \frac{1}{n}\sum_{j=1}^n x_{j,j}$
is the (unique) normalised trace on $\operatorname{M}_n(\mathbf C)$.
\end{Exam}

\begin{Exam}[\textbf{abelian von Neumann algebras}]
\label{abvna}
Let $Z$ be a locally compact space and $\nu$ a positive Radon measure on $Z$.
The space $L^{\infty}(Z,\nu)$ of complex-valued functions on $Z$
which are measurable and $\nu$-essentially bounded 
(modulo equality locally $\nu$-almost everywhere)
is an abelian von Neumann algebra.
The involution is given by $f^*(z) = \overline{f(z)}$.
Any abelian von Neumann algebra is of this form.
\par
If $\nu$ is a probability measure, 
the linear form $\tau_\nu : f \longmapsto \int_Z f(z) d\nu(z)$
is a  trace on $L^{\infty}(Z,\nu)$, normalised in the sense $\tau_\nu(1) = 1$. 
\end{Exam}

\begin{Exam}[\textbf{group von Neumann algebra}]
\label{gvna}
Let $\Gamma$ be a group.
The Hilbert space $\ell^2(\Gamma)$ 
has a scalar product, denoted by  $\langle \cdot \vert \cdot \rangle$,
and a canonical orthonormal basis 
$\left( \delta_\gamma \right)_{\gamma \in \Gamma}$,
where $\delta_\gamma(x)$ is $1$ if $x=\gamma$ and $0$ otherwise.
The \emph{left-regular representation} $\lambda$ of $\Gamma$
on $\ell^2(\Gamma)$ is defined by
$(\lambda(\gamma)\xi)(x) = \xi(\gamma^{-1}x)$
for all $\gamma, x \in \Gamma$ and $\xi \in \ell^2(\Gamma)$.
\par

The \emph{von Neumann algebra $\mathcal N (\Gamma)$ of $\Gamma$} 
is the weak closure in $\mathcal L (\ell^2 (\Gamma))$
of the set of  $\mathbf C$-linear combinations 
$\sum_{\gamma \in \Gamma}^{\text{finite}}  z_\gamma \lambda(\gamma)$;
it is a finite von Neumann algebra.
The involution is given by 
$(z_\gamma \lambda(\gamma))^* = \overline{z_\gamma}\lambda(\gamma^{-1})$.
There is  a canonical trace, given by
$x \longmapsto \langle x \delta_1 \vert \delta_1 \rangle$,
which extends
$\sum_{\gamma \in \Gamma}^{\text{finite}}  
z_\gamma \lambda(\gamma) \longmapsto z_1$.
\par

Moreover,  $\mathcal N (\Gamma)$ is a factor of type $\operatorname{II}_1$
if and only if $\Gamma$ is \emph{icc}\footnote{A
group is \emph{icc} if it is infinite 
and if all its conjugacy classes distinct from $\{1\}$ are infinite.
}
(this is Lemma 5.3.4 in \cite{MuvN--43}, 
see also \cite[chap.\ III, $\S$~7, no 6]{Dixm--57}).
\end{Exam}

\noindent
\textbf{Remarks.}
(a) In the special case of a finite group, 
$\mathcal N (\Gamma)$ of Example \ref{gvna} 
is a finite sum of matrix algebras as in Example \ref{I_nvna}.
In the special case of an abelian group, 
$\mathcal N (\Gamma)$ of Example \ref{gvna} 
is isomorphic, via Fourier transform,
to the algebra of Example \ref{abvna},
with $Z$ the Pontrjagin dual of $\Gamma$ 
(which is a compact abelian group)
and $\nu$ its normalised Haar measure.
\par

(b) The von Neumann algebra $\mathcal N (\Gamma)$ is ``of type $\operatorname{I}$''
if and only if $\Gamma$ has an abelian subgroup of finite index \cite{Thom--64}.
It is ``of type $\operatorname{II}_1$''
if and only if either\footnote{We
denote by $\Gamma_f$ the union 
of the finite conjugacy classes of a group $\Gamma$.
It is easy to check that $\Gamma_f$ is a subgroup,
and it is then obvious that it is a normal subgroup.
}
$[\Gamma : \Gamma_f] = \infty$,
or $[\Gamma : \Gamma_f]  < \infty$ 
and $\vert D\Gamma_f \vert = \infty$ \cite{Kani--69}. 
There exist groups $\Gamma$ such that $\mathcal N (\Gamma)$
is a non-trivial direct product of two two-sided ideals,
one of type $\operatorname{I}$ and the other of type $\operatorname{II}_1$;
see \cite[Theorem 2]{Kapl--51} for the result, 
and \cite{Newm--60} for explicit examples.
\par

(c) Suppose in particular that $\Gamma$ is finitely generated.
If $\Gamma_f$ is of finite index in  $\Gamma$,
then $\Gamma_f$ is also finitely generated and it follows that
$\Gamma$ has an abelian subgroup of finite index.
Thus $\mathcal N (\Gamma)$ is either of type $\operatorname{I}$
(if and only if $\Gamma$ has a free abelian group of finite index)
or of type $\operatorname{II}_1$ 
(if and only if $[\Gamma : \Gamma_f]  = \infty$); 
see \cite{Kani--70}.
\par

(d) Other properties of $\mathcal N (\Gamma)$ are reviewed in \cite{Harp--95}.

\medskip

Let us now recall,  as announced,
some general facts and some terminology.
\begin{itemize}
\item[(i)]
A von Neumann algebra $\mathcal N$ 
inherits \emph{several natural topologies}
from its representations by operators on Hilbert spaces,
including the ``ultraweak topology'' 
(with respect to which the basic examples are separable)
and the ``operator topology'' 
(with respect to which $\mathcal N$ is separable
if and only if it is finite-dimen\-sional).
\item[(ii)]
There is available a \emph{functional calculus}, justified by the \emph{spectral theorem}:
$f(x)$ is well-defined and satisfies natural properties,
for $x \in \mathcal N$ normal ($x^*x = xx^*$)
and $f$ an essentially bounded complex-valued measurable function
on the \emph{spectrum}
\begin{equation*}
\operatorname{sp}(x) \, := \, \{z \in \mathbf C \hskip.1cm \vert \hskip.1cm 
z-x \hskip.2cm \text{is not invertible}\}
\end{equation*}
of $x$. 
More precisely, at least when $\mathcal N$ acts on a separable Hilbert space,
we have for $x$ normal in $\mathcal N$
 a positive regular Borel measure of full support $\nu$ on $\operatorname{sp}(x)$,
 and a natural injective morphism
$L^\infty(\operatorname{sp}(x), \nu) \ni f \longmapsto f(x) \in \mathcal N$
of von Neumann algebras.
\end{itemize}
A von Neumann algebra $\mathcal N$ is \emph{finite}
if, for $x,y \in \mathcal N$, the relation $xy = 1$ implies $yx = 1$.
A \emph{projection} in a von Neumann algebra is a self-adjoint idempotent,
in equations $e = e^* = e^2$.
A  von Neumann algebra $\mathcal N$ if \emph{of type $\operatorname{I}$}
if, for any projection $0 \ne e  \in \mathcal N$,  
there exists a projection $f \in \mathcal N$, $f \ne 0$
such that $fe = ef = f$ and $f \mathcal N f$ is abelian.
A finite von Neumann algebra $\mathcal N$ if \emph{of type $\operatorname{II}_1$}
if, for any projection $0 \ne e \in \mathcal N$,
the subset $e \mathcal N e$ is \emph{not} abelian.
It is known that any finite von Neumann algebra is the direct product
of a finite algebra of type $\operatorname{I}$ 
and an algebra of type $\operatorname{II}_1$.
\begin{itemize}
\item[(iii)]
A \emph{finite trace} on a von Neumann algebra $\mathcal N$
is a linear functional ${\tau : \mathcal N \longrightarrow \mathbf C}$
which is continuous with respect to all the standard topologies on $\mathcal N$
(= which is \emph{normal}, in the standard jargon),
and which satisfies
\subitem(iii$_a$)
$\tau(x^*) = \overline{\tau(x)}$ for all $x \in \mathcal N$,
\subitem(iii$_b$)
$\tau(x^*x) \ge 0$,  for all $x \in \mathcal N$,
\subitem(iii$_c$)
$\tau(xy) = \tau(yx)$ for all $x,y \in \mathcal N$.
\item[]
A trace is \emph{faithful} if $\tau(x^*x) > 0$ whenever $x \ne 0$.
\end{itemize}
It is known that, on a finite von Neumann algebra 
which can be represented on a separable Hilbert space,
there exists a faithful finite normal trace.
Also, any linear form on $\mathcal N$ which is ultraweakly continuous 
and satisfies $(iii)_c$ can be written canonically
as a linear combination of four linear forms satisfying
the three conditions of (iii); this is a \emph{Jordan decomposition}
result of \cite{Grot--57}.
\par

As we will not consider other kind of traces, we use ``trace'' for ``finite trace'' below.

\par

A \emph{factor} is a von Neumann algebra
of which the centre coincides with the scalar multiples of the identity.
A \emph{factor of type $\operatorname{II}_1$} is an infinite dimensional finite factor;
the discovery of such factors is one of the main results of Murray and von Neumann.
\begin{itemize}
\item[(iv)]
Let $\mathcal N$ be a factor of type $\operatorname{II}_1$;
\subitem$\circ$
$\mathcal N$ is a simple ring\footnote{See
\cite[chap.\ III $\S$~5, no 2]{Dixm--57}.
Words are often reluctant to migrate
from one mathematical domain to another.
Otherwise, one could define a factor of type $\operatorname{II}_1$
as an infinite dimensional von Neumann algebra which is central simple.
In the same vein, one could say that von Neumann algebras 
are topologically principal rings;
more precisely, in a von Neumann algebra $\mathcal N$,
any ultraweakly closed left ideal is of the form $\mathcal N e$, 
where $e \in \mathcal N$ is a projection 
(this is a corollary of the von Neumann density theorem
\cite[chap.\ I $\S$~3, no 4]{Dixm--57}).},
\subitem$\circ$
there is a unique normalised\footnote{The
normalisation is most often by $\tau(1) = 1$.
It can be otherwise, for example $\tau(1_n) = n$
on a factor of the form $\operatorname{M}_n(\mathcal N)$,
for some factor $\mathcal N$.} 
trace $\tau$, which is faithful.
\end{itemize}
Thus, on a factor $\mathcal N$ of type $\operatorname{II}_1$,
it is a standard result that
there exists a unique normalised normal trace $\tau$ (in the sense of (iii) above);
but unicity holds in a stronger sense, because
any element in the kernel of $\tau$
is a finite sum of commutators \cite{FaHa--80}.
\par

For a projection $e$, 
the number $\tau(e)$ is called the \emph{von Neumann dimension} of $e$,
or of the Hilbert space $e(\mathcal H)$,
when $\mathcal N$ is understood to be inside some $\mathcal L(\mathcal H)$.

\section{\textbf{The Fuglede-Kadison determinant for finite von Neumann algebras}}
\label{FK1}

In 1952, Fuglede and Kadison defined their \emph{determinant}
\begin{equation}
\label{eqdetgrpeFK}
{\det}^{FK}_{\tau} \, : \, \left\{
\aligned
\operatorname{GL}_1(\mathcal N) \,  &\longrightarrow \hskip2cm \mathbf R_+^*
\\
x \hskip.5cm &\longmapsto  \,
\exp \Big( \tau \big(  \log \big( (x^*x)^{\frac{1}{2}} \big) \big) \Big) 
\endaligned \right. 
\end{equation}
which is a partial analogue of (\ref{eqdetgrpeclassique}).
The number ${\det}^{FK}_{\tau} (x)$ is well-defined by functional calculus,
and most of the work in \cite{FuKa--52} is for showing that ${\det}^{FK}_{\tau}$
is a homomorphism of groups.
For the definition given below in Section \ref{FK2}, 
it will be the opposite:
some work to show that the definition makes sense,
but a very short proof to show it defines a group homomorphism.

In the original paper, $\mathcal N$ is a factor of type $\operatorname{II}_1$,
and $\tau$ is its unique trace with $\tau(1) = 1$;
but everything carries over to the case of a von Neumann algebra
and a normalised trace \cite[chap.\ I, $\S$~6, no 11]{Dixm--57}.
Besides being a group homomorphism, ${\det}^{FK}_{\tau}$ has
the following properties:
\begin{equation*}
\aligned
\circ \hskip.2cm  
& {\det}^{FK}_{\tau}(e^y) \, = \, \vert e^{\tau (y)} \vert
\, = \, e^{\operatorname{Re}(\tau(y))}
\hskip.2cm \text{for all} \hskip.2cm
y \in \mathcal N
\\
& \text{and in particular} \hskip.2cm 
{\det}^{FK}_{\tau}(\lambda 1) \, = \, \vert \lambda \vert
\hskip.2cm \text{for all} \hskip.2cm 
\lambda \in \mathbf C ,
\\
\circ \hskip.2cm  
& {\det}^{FK}_{\tau} (x) = {\det}^{FK}_{\tau} \big( (x^*x)^{\frac{1}{2}} \big)
\hskip.5cm \text{for all} \hskip.2cm 
x \in \operatorname{GL}_1(\mathcal N) 
\\
& \text{and in particular,} \hskip.2cm {\det}^{FK}_{\tau} (x) = 1
\hskip.2cm \text{for all} \hskip.2cm
x \in \operatorname{U}_1(\mathcal N) .
\endaligned
\end{equation*}
For a $*$-ring $\mathcal R$ with unit,  
$\operatorname{U}_1(\mathcal R)$ denotes
its \emph{unitary group}, defined to be
$\{x \in \mathcal \mathcal R \hskip.1cm \vert \hskip.1cm x^*x = xx^* = 1 \}$.
\par

Instead of (\ref{eqdetgrpeFK}), we could equally view ${\det}^{FK}_{\tau}$
as a family of homomorphisms
$\operatorname{GL}_n(\mathcal N) \longrightarrow  \mathbf R_+^*$,
one for each $n \ge 1$;
if the traces on the $\operatorname{M}_n(\mathcal N)$ 's
are normalised by $\tau(1_n) = n$, we have
${\det}^{FK}_{\tau}(\lambda 1_n) = \vert \lambda \vert^n$.
More generally, for any projection 
$e \in \operatorname{M}_n(\mathcal N)$,
we have a von Neumann algebra 
$\operatorname{M}_e(\mathcal N) := e\operatorname{M}_n(\mathcal N)e$,
and a Fuglede-Kadison determinant
${\det}^{FK}_{\tau} : \operatorname{GL}_e(\mathcal N) \longrightarrow \mathbf R_+^*$
defined on its group of units.
\par

There are extensions of ${\det}^{FK}_{\tau}$ to non-invertible elements,
but this raises some problems, and technical difficulties.
Two extensions are discussed in \cite{FuKa--52}:
the ``algebraic extension'' for which 
the determinant vanishes on singular elements
(this is not mentionned again below),
and the ``analytic extension''
which relies on Formula (\ref{eqdetgrpeFK}),
in which one should understand
\begin{equation}
\label{analyticextension}
{\det}^{FK}_{\tau}(x) = 
\exp \Big( \tau \big( \log ( (x^*x)^{\frac{1}{2}} ) \big) \Big)  =
\exp
\int_{\operatorname{sp}((x^*x)^{1/2})} 
\ln \lambda \hskip.1cm d \tau(E_\lambda) ,
\end{equation}
where $(E_\lambda)_{\lambda \in \operatorname{sp}((x^*x)^{1/2})}$
holds for the spectral measure of $(x^*x)^{1/2}$;
of course $\exp ( - \infty ) = 0$.
(Note that we write ``$\log$'' for logarithms of matrices and operators,
and ``$\ln$'' for logarithms of numbers.)
For example, if $x$ is such that 
there exists a projection $e$ with $x = x(1-e)$ and $\tau(e) > 0$,
we have ${\det}^{FK}_{\tau}(x) = 0$.
For all $x,y \in \mathcal N$, we have
\begin{equation*}
\aligned
{\det}^{FK}_{\tau}((x^*x)^{1/2}) \, &= \,  
\lim_{\epsilon \to 0+} {\det}^{FK}_{\tau}((x^*x)^{1/2} + \epsilon 1)
\\
{\det}^{FK}_{\tau}(x) {\det}^{FK}_{\tau}(y) \, &= \,  {\det}^{FK}_{\tau}(xy)
\endaligned
\end{equation*}
(see \cite{FuKa--52}, respectively Lemma 5 and Page 529).
But an element $x$ with ${\det}^{FK}_{\tau}(x) \ne 0$ need not be invertible,
and no extension $\mathcal N \longrightarrow \mathbf R_+$
of the mapping ${\det}^{FK}_{\tau}$ of (\ref{eqdetgrpeFK})
is norm-continuous \cite[Theorem~6]{FuKa--52}.

We will discuss another extension ${\det}^{FKL}_{\tau}$
to singular elements, in Section \ref{baratinL2torsion}.

\medskip

More generally, ${\det}^{FK}_{\tau}(x)$ can be defined for $x$
an operator ``affiliated'' to $\mathcal N$, 
and also for traces which are \emph{semi-finite}
rather than finite as above.
See \cite{Grot--55, Arve--67, Fac--82b, Fack--83, Brow--86, HaSc--09},
among others.
We will not have any further comment on this part of the theory.

\begin{Exam}[\textbf{Fuglede-Kadison determinant 
for $\operatorname{M}_n(\mathbf C)$}]
Let $\mathcal N = \operatorname{M}_n(\mathbf C)$
be the factor of type $\operatorname{I}_n$,
as in Example \ref{I_nvna},
let $\det$ be the \emph{usual} determinant, 
and let  $\tau : x \longmapsto \frac{1}{n}\sum_{j=1}^n x_{j,j}$
be the  trace normalised by $\tau(1_n) = 1$.
Then
%
%
\begin{equation}
{\det}^{FK}_{\tau}(x) \, = \,  \vert \det (x) \vert^{1/n} 
\, = \,  \big( \det ( (x^*x)^{1/2}) \big)^{1/n}
\end{equation}
for all $x \in \operatorname{M}_n(\mathbf C)$.
\end{Exam}

\begin{Exam}[\textbf{Fuglede-Kadison determinant for abelian von Neumann algebras}]
Let $L^{\infty}(Z,\nu)$ and $\tau_\nu$ be as in Example \ref{abvna},
with $\nu$ a probability measure.
The corresponding Fuglede-Kadison determinant is given by
\begin{equation}
\label{!!}
{\det}^{FK}_{\tau} (f) \, = \,  \exp \int_Z \ln \vert f (z) \vert d \mu (z)
\, \in \,  \mathbf R_+ .
\end{equation}
In (\ref{!!}), observe  that $\ln \vert f (z) \vert$ is bounded above on $Z$,
because $\vert f(z) \vert \le \Vert f \Vert_{\infty} < \infty$
for $\nu$-almost all $z$.
However $\vert f(z) \vert$  need not be bounded away from $0$,
so that $\ln \vert f(z) \vert = -\infty$ occurs.
If the value of the integral is $-\infty$, then
${\det}^{FK}_{\tau} (f) = \exp(-\infty) = 0$.
\end{Exam}

\medskip

Consider an integer $d \ge 1$ and the von Neumann algebra 
$\mathcal N (\mathbf Z^d)$ of the free abelian group of rank $d$.
Fourier transform provides an isomorphism of von Neumann algebras
\begin{equation*}
\mathcal N (\mathbf Z^d) 
\overset{\approx}{\longrightarrow}
\operatorname{L}^{\infty}(T^d, \nu),
\hskip.5cm x \longmapsto \hat x ,
\end{equation*}
where $\nu$ denotes the normalised Haar measure
on the $d$-dimensional torus $T^d$.
Moreover, the composition of this isomorphism
with the trace $\tau_\nu$ of Example \ref{abvna}  
is the canonical trace on $\mathcal N (\mathbf Z^d)$,
in the sense of Example \ref{gvna}.

\begin{Exam}[\textbf{Fuglede-Kadison determinant and Mahler measure}]
Let $x$ be a finite linear combination
$\sum_{n \in \mathbf Z^d}^{finite} z_n \lambda(n) \in \mathcal N (\mathbf Z^d)$,
so that $\hat x \in L^{\infty}(T^d,\nu)$ is a trigonometric polynomial.
Then the $\tau_\nu$-Fuglede-Kadison determinant of $x$ 
is given by the \emph{exponential Mahler measure} of $\hat x$:
\begin{equation*}
{\det}^{FK}_{\tau_\nu} (x) \, = \, M(\hat x) \, := \,  
\exp \int_{T^d} \ln \vert \hat x (z) \vert d \nu (z) .
\end{equation*}
\end{Exam}
In the one-dimensional case ($d=1$), if
\begin{equation*}
\hat x(z) = a_0 + a_1z + \cdots + a_s z^s = a_s \prod_{j=1}^s (z-\xi_j) ,
\hskip.5cm \text{with} \hskip.2cm a_0 a_s \ne 0 ,
\end{equation*} 
a computation shows that
\begin{equation*}
\int_{T} \ln \vert \hat x (z) \vert d \nu (z) = 
\int_0^1  \ln \vert \hat x (e^{2 \pi i \alpha}) \vert d\alpha =
\ln  \vert a_s \vert + \sum_{j=1}^s \max \{1, \vert \xi_j \vert \}
\end{equation*}
(see \cite[Proposition 16.1]{Schm--95} or \cite[Pages 135--7]{Luck--02}).
\par

Mahler measures occur in particular as entropies of $\mathbf Z^d$-actions 
by automorphisms of compact groups.
More precisely, for $x \in \mathbf Z [\mathbf Z^d]$,
which can be viewed as the inverse Fourier transform of a trigonometric polynomial,
the group $\mathbf Z^d$
acts naturally on the quotient  $\mathbf Z [\mathbf Z^d]/(x)$
of the group ring by the principal ideal $(x)$,
hence on the Pontryagin dual $\left( \mathbf Z [\mathbf Z^d] / (x) \right)^{\widehat{}}$
of this countable abelian group, which is a compact abelian group. 
For example, if $x(z) = 1 + z - z^2 \in \mathbf Z [z, z^{-1}] \approx \mathbf Z [\mathbf Z]$,
then $\left( \mathbf Z [\mathbf Z] / (x) \right)^{\widehat{}} \approx T^2$,
and the corresponding action of the generator of $\mathbf Z$ on $T^2$
is described by the matrix 
$\begin{pmatrix}
0 & 1 \\ 1 & 1
\end{pmatrix}$ \cite[Example 5.3]{Schm--95}.
Every action of $\mathbf Z^d$ by automorphisms of a compact abelian group
arises as above from some $x \in \mathbf Z [\mathbf Z^d]$. 
More on this in \cite{LiSW--90, Schm--95, Deni--06}.
\par

The logarithm of the Fuglede-Kadison determinant occurs also 
in the definition of a  ``tree entropy'', namely in the asymptotics
of the number of spanning trees in large graphs
\cite{Lyon--05, Lyon--10}.

\section{\textbf{A motivating question}}
\label{motivating}

\emph{
It is natural  to ask why $\mathbf R_+^*$ appears on
the right-hand side of  (\ref{eqdetgrpeFK}),
even  though $\mathcal N$ is a \emph{complex} algebra,
for example a $\operatorname{II}_1$-factor,
whereas $\mathbf C^*$ appears on the right-hand side of (\ref{detexpC2})
when $\mathcal N = \operatorname{M}_n(\mathbf C)$.
}

\medskip

This is not due to some shortsightedness of Fuglede and Kadison.
Indeed, for $\mathcal N$ a factor of type $\operatorname{II}_1$,
it has been shown
that the Fuglede-Kadison determinant provides an \emph{isomorphism}
from the abelianised group
$\operatorname{GL}_1(\mathcal N) / D\operatorname{GL}_1(\mathcal N)$
onto $\mathbf R_+^*$.
In other words:

\begin{Prop}[\textbf{Properties of operators with trivial Fuglede-Ka\-dison determinant
in a factor of type $\operatorname{II}_1$}]
\label{SL}
Let $\mathcal N$ be a factor of type $\operatorname{II}_1$.
\begin{itemize}
\item[(i)]
Any element in $\operatorname{U}_1(\mathcal N)$ 
is a product of finitely many multiplicative commutators
of unitary elements.
\item[(ii)]
The kernel $\operatorname{SL}_1(\mathcal N)$
of the homomorphism (\ref{eqdetgrpeFK})
coincides with the group of commutators 
of $\operatorname{GL}_1(\mathcal N)$.
\end{itemize}
\end{Prop}

Property (i) is due to Broise \cite{Broi--67}.
It is moreover known that any proper normal subgroup of 
$\operatorname{U}_1(\mathcal N)$
is contained in its center, which is
$\{ \lambda \operatorname{id} \hskip.1cm \vert \hskip.1cm 
\lambda \in \mathbf C^*, \vert \lambda \vert = 1 \} \approx \mathbf R / \mathbf Z$
\cite[Proposition 3 and its proof]{Harp--79};
this sharpens an earlier result 
on the classification of norm-closed normal subgroups of 
$\operatorname{U}_1(\mathcal N)$ \cite[Theorem 2]{Kadi--52}.
\par
Property (ii) is  \cite[Proposition 2.5]{FaHa--80}.
It follows that the quotient of $\operatorname{SL}_1(\mathcal N)$ by its center
(which is the same as the centre of $\operatorname{U}_1(\mathcal N)$)
is simple, as an abstract group \cite[Corollary 6.6, Page 123]{Lans--70}.

\medskip

As a kind of answer to our motivating question,
we will see below that, 
when the Fuglede-Kadison definition 
is adapted to a \emph{separable} Banach algebra,
the right-hand side of the homomorphism analogous to (\ref{eqdetgrpeFK})
is necessarily a quotient of the additive group $\mathbf C$
by a \emph{countable} subgroup.
For example, when $A = \operatorname{M}_n(\mathbf C)$, this quotient is
$\mathbf C / 2i\pi \mathbf Z \overset{\exp(\cdot)}{\approx} \mathbf C^*$,
see Corollary \ref{classicalA}.
On the contrary, when $A$ is a $\operatorname{II}_1$-factor 
(not separable as a Banach algebra), this quotient is
$\mathbf C / 2i\pi \mathbf R 
\overset{\exp(\operatorname{Re}(\cdot))}{\approx} \mathbf R_ +^*$,
see Corollary \ref{A=M}.
The case of a separable Banach algebra can sometimes be seen
as providing an interpolation between the two previous cases,
see Remark \ref{interp}.

\section{\textbf{A reminder on $K_0$, $K_1$, $K_1^{\operatorname{top}}$,
and Bott periodicity}}
\label{K}

\subsection{On $K_0(\mathcal R)$ and $K_0(A)$}
Let $\mathcal R$ be a ring, say with unit to simplify several small technical points. 
Let us first recall one definition of the abelian group $K_0(\mathcal R)$
of K-theory.
\par

We have a nested sequence of rings of matrices and
(non-unital) ring homomorphisms
\begin{equation}
\aligned
\mathcal R = 
\operatorname{M}_1(\mathcal R) \subset \cdots \subset 
\operatorname{M}_n(\mathcal R) 
&\subset \operatorname{M}_{n+1}(\mathcal R) \subset \cdots
\\
&\subset \operatorname{M}_\infty(\mathcal R) 
:= \bigcup_{n \ge 1}\operatorname{M}_n(\mathcal R),
\endaligned
\end{equation}
where the inclusions at finite stages are given by
$x \longmapsto 
\begin{pmatrix}
x & 0 \\ 0 & 0
\end{pmatrix}$.
\par

An \emph{idempotent} in $\operatorname{M}_\infty(\mathcal R)$ 
is an element $e$ such that $e^2 = e$.
Two idempotents $e, f \in \operatorname{M}_\infty(\mathcal R)$ are \emph{equivalent}
if there exist $n \ge 1$ and $u \in \operatorname{GL}_n(\mathcal R)$
such that $e,f \in \operatorname{M}_n(\mathcal R)$ and $f = u^{-1}eu$.
Define an \emph{addition} on equivalence classes of idempotents, by
\begin{equation}
\label{plusinK_0}
\aligned
(\text{class of } e \in \operatorname{M}_k(\mathcal R)) 
&+ (\text{class of } f \in  \operatorname{M}_\ell(\mathcal R))
\\
& \hskip1cm
\, = \, 
\text{class of } 
e \oplus f
\in \operatorname{M}_{k + \ell}(\mathcal R)
\endaligned
\end{equation}
where $e \oplus f$ denotes the matrix
$\begin{pmatrix}
e & 0 \\ 0 & f
\end{pmatrix}$.
Two idempotents $e,f \in \operatorname{M}_\infty(\mathcal R)$ 
are \emph{stably equivalent}
if there exists an idempotent $g$ 
such that the classes of $e \oplus g$ and $f \oplus g$ are equivalent;
we denote by $[e]$ the stable equivalence class of an idempotent $e$.
The set of stable equivalence classes of idempotents,
with the addition defined by $[e] + [f] := [e \oplus f]$, 
is a semi-group.
The Grothendieck group $K_0(\mathcal R)$ of this semi-group 
is the set of formal differences
$[e]-[e']$, up to the equivalence defined by $[e]-[e'] \sim [f]-[f']$
if $[e] + [f'] = [e'] + [f]$.

Note that $K_0$ is a functor: 
to any (unital) ring homomorphism $\mathcal R \longrightarrow \mathcal R'$ 
corresponds a natural homomorphism 
$K_0(\mathcal R) \longrightarrow K_0(\mathcal R')$
of abelian groups.
Note also the isomorphism 
$K_0(\operatorname{M}_n(\mathcal R)) \approx K_0(\mathcal R)$,
which is a straightforward consequence of the definition and of the isomorphisms
$\operatorname{M}_k(\operatorname{M}_n(\mathcal R)) \approx
\operatorname{M}_{kn}(\mathcal R)$.

(To an idempotent $e \in \operatorname{M}_\infty(\mathcal R)$ is associated 
a $\mathcal R$-linear mapping $\mathcal R^n \to \mathcal R^n$ for $n$ large enough,
of which the image 
is a projective $\mathcal R$-module of finite rank.
From this it can be checked that the definition of $K_0(\mathcal R)$ given above 
coincides with another standard definition, 
in terms of projective modules of finite rank.
Details in \cite[Chap.\ 1]{Rose--94}.)

\medskip

Rather than a general ring $\mathcal R$, consider now the case
of a complex Banach algebra $A$ with unit.
For each $n \ge 1$, the matrix algebra $\operatorname{M}_n(A)$ 
is again a Banach algebra, for some appropriate norm,
and we can furnish $\operatorname{M}_\infty(A)$
with the inductive limit topology.
The following is rather easy to check,
see e.g. \cite[Pages 25--27]{Blac--86}:
two idempotents $e,f \in \operatorname{M}_\infty(A)$ 
are equivalent if and only if  there exists a continuous path
\begin{equation*}
[0,1] \longrightarrow \{\text{idempotents of } \operatorname{M}_\infty(A) \},
\hskip.5cm \alpha \longmapsto e_\alpha
\end{equation*}
such that $e_0 = e$ and $e_1 = f$.
This has the following consequence:

\begin{Prop}
\label{Asep}
If the Banach algebra $A$ is separable, 
the abelian group $K_0(A)$ is countable.
\end{Prop}

\begin{Prop}
\label{Kfacteurfini}
If $\mathcal N$ is a factor of type $\operatorname{II}_1$, 
then $K_0(\mathcal N) \approx \mathbf R$ is uncountable.
\par
Indeed, if $\tau$ denotes the canonical trace on $\mathcal N$, 
the mapping which associates to 
the class of a self-adjoint idempotent $e$ in $\mathcal N$
its von Neumann dimension $\tau(e) \in [0,1]$ extends to an isomorphism
${K_0(\mathcal N) \overset{\approx}{\longrightarrow} \mathbf R}$.
\end{Prop}

\noindent
\emph{On the proof~:} this follows from the ``comparison of projections''  
in von Neumann algebras \cite[chap.\ III, $\S$~2, no 7]{Dixm--57}.
\hfill $\square$

\medskip

For historical indications on the early connections between K-theory
and operator algebras, which goes back to the mid 60's, see \cite{Rose--05}.

\subsection{On $K_1(\mathcal R)$}
\label{K1alg}
For any ring $\mathcal R$ with unit,
we have a nested sequence of group homomorphisms
\begin{equation}
\label{eqGLpourK1}
\aligned
\mathcal R^* =
\operatorname{GL}_1(\mathcal R) \subset \cdots \subset 
\operatorname{GL}_n(\mathcal R) &\subset 
\operatorname{GL}_{n+1}(\mathcal R) \subset \cdots
\\
&\subset \operatorname{GL}_\infty(\mathcal R) 
:= \bigcup_{n \ge 1} \operatorname{GL}_n(\mathcal R),
\endaligned
\end{equation}
where the inclusions at finite stages are given by
$x \longmapsto 
\begin{pmatrix}
x & 0 \\ 0 & 1
\end{pmatrix}$.

By definition, 
\begin{equation}
K_1(\mathcal R) \, = \,
\operatorname{GL}_\infty(\mathcal R) / D\operatorname{GL}_\infty(\mathcal R) 
\end{equation}
is an abelian group, usually written additively. 
Note that $K_1$ is a functor from rings to abelian groups.
\par

For a commutative ring $\mathcal R$,
the classical determinant provides a homomorphism
$K_1(\mathcal R) \longrightarrow \mathcal R^*$;
it is an isomorphism in several important cases,
for example when $\mathcal R$ a field,
or the ring of integers in a finite extension of $\mathbf Q$
\cite[$\S$~3]{Miln--71}.
In general ($\mathcal R$ commutative or not), 
the association of an element in $K_1(\mathcal R)$
to a matrix in $\operatorname{GL}_\infty(\mathcal R)$ 
can be viewed as a kind of determinant, or rather of a log of a determinant
since $K_1(\mathcal R)$ is written additively.
Accordingly,  the torsion defined in (\ref{deftauMilnor}) below
can be viewed as an alternating sum of log of determinants;
we will remember this when defining the $L^2$-torsion 
in Equation (\ref{defL2torsion}).
\par

Let $\mathcal R$ be again an arbitrary ring with unit.
The  \emph{reduced} $K_1$-group is the quotient $\overline{K}_1(\mathcal R)$
of $K_1(\mathcal R)$ by the image of the natural homomorphism
$\{1, -1\} \subset \operatorname{GL}_1(\mathcal R) 
\subset \operatorname{GL}_\infty(\mathcal R) \longrightarrow K_1(\mathcal R)$.
\par

In case $\mathcal R = \mathbf Z [\Gamma]$ 
is the integral group ring of a group $\Gamma$, 
the \emph{Whitehead group} $\operatorname{Wh}(\Gamma)$
is the cokernel $K_1(\mathbf Z [\Gamma]) / \langle \pm 1, \Gamma \rangle$
of the natural homomorphism
$\Gamma \subset \operatorname{GL}_1(\mathbf Z [\Gamma])
\longrightarrow K_1(\mathbf Z [\Gamma])
\longrightarrow \overline{K}_1(\mathbf Z [\Gamma])$.
\par

When $\Gamma$ is finitely presented,
there is a different (but equivalent) definition of $\operatorname{Wh}(\Gamma)$, 
with geometric content.
In short, let $L$ be a connected finite CW-complex with $\pi_1(L) = \Gamma$.
One defines a group $\operatorname{Wh}(L)$ of appropriate equivalence classes
of pairs $(K,L)$, with $K$ a finite CW-complex containing $L$
in such a way that the inclusion $L \subset K$
is a homotopy equivalence.
The unit is represented by pairs $L \subset K$
for which the inclusion is a \emph{simple} homotopy equivalence.
It can be shown that the functors $L \longrightarrow \operatorname{Wh}(L)$
and $L \longrightarrow \operatorname{Wh}(\pi_1(L))$
are naturally equivalent \cite[$\S$~6 and Theorem 21.1]{Cohe--73}.

Examples:  
$\operatorname{Wh}(\mathbf Z^d) = 0$ 
for free abelian groups $\mathbf Z^d$
and $\operatorname{Wh}(F_d) = 0$ for 
free groups $F_d$.
For finite cyclic groups:
$\operatorname{Wh}(\mathbf Z / q \mathbf Z)$ is a free abelian group
of finite rank for all $q \ge 1$, 
and is the group $\{0\}$ if and only if $q \in \{1,2,3,4,6\}$.

\medskip  

From the standard references, let us quote
\cite{RhMK--67}, \cite{Miln--66},
\cite{Miln--71},
\cite{Cohe--73}, and \cite{Tura--01}.

\subsection{On  $K_1^{\operatorname{top}}(A)$,
and on $K_0(A)$ viewed as a fundamental group}
\label{K0=pi1}
Let $A$ be a Banach algebra with unit.
For each $n \ge 1$, the group $\operatorname{GL}_n(A)$ 
is an open subset of the Banach space $\operatorname{M}_n(A)$, 
and the induced topology makes it a topological group.
The group $\operatorname{GL}_\infty(A)$ of (\ref{eqGLpourK1})
is also a topological group, for the inductive limit topology;
we denote by $\operatorname{GL}_\infty^0(A)$ its connected component.
\par

It is a simple consequence of the classical ``Whitehead lemma''
that, for any Banach algebra, the group  $D\operatorname{GL}_\infty(A)$ 
is perfect and coincides with $D\operatorname{GL}_\infty^0(A)$;
see for example \cite[Appendix]{HaS--85}.
In particular, $D\operatorname{GL}_\infty(A) \subset
\operatorname{GL}_\infty^0(A)$,
to that the quotient group
\begin{equation}
K_1^{\operatorname{top}}(A) \, := \,  \pi_0\left( \operatorname{GL}_\infty(A) \right)
\, = \, \operatorname{GL}_\infty(A) /  \operatorname{GL}_\infty^0(A)
\end{equation}
is commutative.
Note that $ \operatorname{GL}_1(A) /  \operatorname{GL}_1^0(A)$
need not be commutative \cite{Yuen--73},
even if its image in $\operatorname{GL}_\infty(A) /  \operatorname{GL}_\infty^0(A)$
is always commutative.
\par

Moreover,
we have a natural quotient homomorphism
\begin{equation}
\label{K1etK1top}
\operatorname{GL}_\infty(A) / D\operatorname{GL}_\infty^0(A)
= K_1(A) 
 \longrightarrow 
K_1^{\operatorname{top}}(A) = 
\operatorname{GL}_\infty(A) / \operatorname{GL}_\infty^0(A)
\end{equation}
which is onto.
It is an isomorphism if and only if the group $\operatorname{GL}_\infty^0(A)$
is perfect; this is the case if $A$ is an infinite simple C$^*$-algebra,
for example if $A$ is one of the Cuntz algebras $O_n$ briefly mentionned below.
\par

If the Banach algebra $A$ is separable, the group $K_1^{\operatorname{top}}(A)$
is countable (compare with Proposition \ref{Asep}).
\par

To an idempotent $e \in \operatorname{M}_n(A)$, we can associate the loop
\begin{equation}
\label{xie}
\xi_e \, : \, 
\left\{
\aligned 
 {} 
 [0,1] \hskip.2cm &\longrightarrow \hskip2cm 
\operatorname{GL}_n(A)  \subset \operatorname{GL}_\infty(A)
\\ 
\alpha \hskip.5cm &\longmapsto \hskip.2cm \exp (2 \pi i \alpha e) \, = \, 
\exp (2\pi i \alpha)e + (1-e) ;
\endaligned
\right.
\end{equation}
note that $\xi_e(0) = \xi_e(1) = 1$.
If two idempotents $e$ and $f$ have the same image in $K_0(A)$,
it is easy to check that 
$\xi_e$ and $\xi_f$ are homotopic loops.
It is a fundamental fact, which is a form of \emph{Bott periodicity},
that the assignment $e \longmapsto \xi_e$
extends to a group isomorphism
\begin{equation}
\label{bott}
K_0(A) \overset{\approx}{\longrightarrow} 
\pi_1 \left( \operatorname{GL}_\infty^0(A) \right) ;
\end{equation}
see \cite[Theorem III.1.11]{Karo--78} or \cite[Chapter 9]{Blac--86}.
The terminology is due to a generalisation of (\ref{bott}):
$K_i^{\operatorname{top}}(A) \approx K_{i+2}^{\operatorname{top}}(A)$
for any integer $i \ge 0$;
by definition, 
$K_i^{\operatorname{top}}(A) = \pi_{i-1} \left( \operatorname{GL}_\infty(A) \right)$,
for all $i \ge 1$, and $K_0^{\operatorname{top}}(A) = K_0(A)$.

\subsection{A few standard examples}
\label{afewstandard}
Let $A = \mathcal C (T)$ be the Banach algebra of continuous functions
on a compact space $T$.
Then $K_0(A) = K^0(T)$ and $K_1^{\operatorname{top}}(A) = K^1(T)$,
where $K^0(T)$ and $K^1(T)$ stand for the (Grothendieck)-Atiyah-Hirzebruch-Bott
K-theory groups of the topological space $T$, 
defined in terms of complex vector bundles.
For example, if $T$ is a sphere, we have
\begin{equation*}
\aligned
K_0(\mathcal C (\mathbf S^{2m})) \approx \mathbf Z^2 ,
\hskip.5cm &\hskip.5cm
K_1^{\operatorname{top}}(\mathcal C (\mathbf S^{2m})) = 0 ,
\\
K_0(\mathcal C (\mathbf S^{2m+1})) \approx \mathbf Z ,
\hskip.5cm &\hskip.5cm
K_1^{\operatorname{top}}(\mathcal C (\mathbf S^{2m+1})) \approx \mathbf Z ,
\endaligned
\end{equation*}
for all  $m \ge 0$.
If $T$ is a compact CW-complex without cells of odd dimension,
then $K_1^{\operatorname{top}}(\mathcal C (T)) = 0$.
\par

Let $A$ be an \emph{AF-algebra}, 
namely a C$^*$-algebra which contains a nested sequence 
$A_1 \subset \cdots \subset A_n \subset A_{n+1} \subset \cdots$
of finite-dimensional sub-C$^*$-algebras with $\bigcup_{n \ge 1} A_n$ dense in $A$.
Then $K_0(A)$ is rather well understood, and
$
K_1^{\operatorname{top}}(A) = 0 .
$
The group $K_0(A)$ is the basic ingredient in Elliott's classification of AF-algebras,
from the 1970's; this was the beginning of a long and rich story,
with a numerous offspring, see 
\cite[Chapter 7]{Blac--86}, \cite{Rord--02}, and \cite{ElTo--08}.
Here is a particular case, the so-called \emph{CAR algebra},
or C$^*$-algebra of the \emph{Canonical Anticommutation Relations}:
it is the C$^*$-closure of the limit of the inductive system of finite matrix algebras
\begin{equation*}
\mathbf C \, \subset \, \cdots \, \subset \, 
\operatorname{M}_{2^n}(\mathbf C)  \, \subset \, 
\operatorname{M}_{2^{n+1}}(\mathbf C)  \, \subset \, 
\cdots 
\end{equation*}
where the inclusions  are given by
$x \longmapsto 
\begin{pmatrix}
x & 0 \\ 0 & x
\end{pmatrix}$.
For this,
\begin{equation*}
K_0(CAR) \, = \,  \mathbf Z [1/2]
\hskip.5cm \text{and} \hskip.5cm
K_1^{\operatorname{top}}(CAR) \, = \, 0 
\end{equation*}
(for $K_1$ of $CAR$ and a few other AF-algebras, 
see Subsection \ref{subsectionsharp}).
\par

The Jiang-Su algebra $\mathcal Z$ is 
a simple infinite-dimensional C$^*$-algebra with unit
which plays an important role in Elliott's classification program of C$^*$-algebras.
It has the same K-theory as $\mathbf C$ \cite{JiSu--99}. 

\par

The \emph{reduced C$^*$-algebra} of a group $\Gamma$
is the \emph{norm}-closure $C^*_\lambda (\Gamma)$
of the algebra 
$\left\{ \sum_{\gamma \in \Gamma}^{\text{finite}} z_\gamma \lambda(\gamma) \right\}$,
see Example \ref{gvna}, in the algebra of all bounded operators
on $\ell^2(\Gamma)$. For the free groups $F_d$ (non-abelian free groups if $d \ge 2$),
we have \cite{PiVo--82}
\begin{equation*}
K_0(C^*_\lambda (F_d)) \,  \approx \, \mathbf Z
\hskip.5cm \text{and} \hskip.5cm
K_1^{\operatorname{top}}(C^*_\lambda (F_d)) \, \approx \,  \mathbf Z^d .
\end{equation*}
\par

For a so-called \emph{irrational rotation C$^*$-algebra} $A_\theta$,
generated by two unitaries $u,v$ satisfying the relation $uv = e^{2 \pi i \theta}vu$,
where $\theta \in [0,1]$ with $\theta \notin \mathbf Q$, 
we have \cite{PiVo--80}
\begin{equation*}
K_0(A_\theta) \,  \approx \, \mathbf Z^2
\hskip.5cm \text{and} \hskip.5cm
K_1^{\operatorname{top}}(A_\theta) \, \approx \,  \mathbf Z^2.
\end{equation*}
\par

For the infinite  \emph{Cuntz algebras} $O_n$,
generated by $n \ge 2$ elements $s_1, \hdots, s_n$ 
satisfying $s_j^*s_k = \delta_{j,k}$ and $\sum_{j=1}^n s_j s_j^* = 1$,
we have \cite{Cuntz--81}
\begin{equation*}
K_0(O_n) \,  \approx \, \mathbf Z / (n-1)\mathbf Z  
\hskip.5cm \text{and} \hskip.5cm
K_1^{\operatorname{top}}(O_n) \, = 0 .
\end{equation*}
\par

For $\mathcal N$ a factor of type $\operatorname{II}_1$, we have
\begin{equation*}
K_0(\mathcal N) \, \approx \,  \mathbf R 
\hskip.5cm  \text{and} \hskip.5cm
K_1(\mathcal N) \, = \,   \mathbf R_+^* .
\end{equation*}
For $K_0$, see Proposition \ref{Kfacteurfini};
for $K_1$, see \cite[already cited for Proposition \ref{SL}.ii]{FaHa--80}. 
More generally, for $\mathcal N$ a von Neumann algebra 
of type  $\operatorname{II}_1$, 
with centre denoted by $\mathcal Z$, we have
\begin{equation*}
K_0(\mathcal N) \, \approx \,  
\{ z \in \mathcal Z \hskip.1cm \vert \hskip.1cm z^*=z \} ,
\end{equation*} 
where the right-hand side is viewed as a group for the addition, and
\begin{equation*}
K_1(\mathcal N) \, \approx \,  
\{ z \in \mathcal Z \hskip.1cm \vert \hskip.1cm
z \ge \epsilon > 0 \}
\hskip1cm \text{($\epsilon$ depends on $z$)},
\end{equation*}
where the right-hand side is viewed as a group for the multiplication;
see \cite{LuRo--93} or \cite[Section 9.2]{Luck--02}.
For any von Neumann algebra $\mathcal N$
\begin{equation*}
K_1^{\operatorname{top}}(\mathcal N) \, = \,  0 ,
\end{equation*}
because $\operatorname{GL}_n(\mathcal N)$ is connected for all $n \ge 1$;
indeed, by polar decomposition and functional calculus, 
any $x \in \operatorname{GL}_n(\mathcal N)$ is of the form $\exp(a)\exp(ib)$,
with $a,b$ self-adjoint in $\operatorname{M}_n(\mathcal N)$,
so that $x$ is connected to $1$
by the path $\alpha \longmapsto \exp(\alpha a)\exp(i \alpha b)$.
\par

\subsection{The topology of the group of units in a factor of type $\operatorname{II}_1$,
and Bott periodicity}
If $\mathcal N$ is a factor of type $\operatorname{II}_1$,
the isomorphism (\ref{bott}) of Bott periodicity shows that
\begin{equation*}
\pi_1 \left( \operatorname{GL}_\infty(\mathcal N) \right) 
\, \approx \, 
K_0(\mathcal N)
\, \approx \ \mathbf R .
\end{equation*}
Thus, by Bott periodicity, 
\begin{equation*}
\pi_{2j} \left( \operatorname{GL}_\infty(\mathcal N) \right)  \, = \, 0
\hskip.5cm \text{and} \hskip.5cm
\pi_{2j+1} \left( \operatorname{GL}_\infty(\mathcal N) \right) \, \approx \, \mathbf R
\end{equation*}
for all $j \ge 0$.
\par

For $\pi_1$, it is known more precisely that
$\pi_1\left( \operatorname{GL}_n(\mathcal N) \right) \approx \mathbf R$,
and that the embedding of $\operatorname{GL}_n(\mathcal N)$ 
into $\operatorname{GL}_{n+1}(\mathcal N)$
induces the identity on $\pi_1$, for all $n \ge 1$
\cite{ArSS--71, Hand--78}.
Note that, still for the norm topology, 
polar decomposition shows that
the unitary group $\operatorname{U}_1(\mathcal N)$
is a deformation retract of $\operatorname{GL}_1(\mathcal N)$;
in particular, we have also
$\pi_1(\operatorname{U}_1(\mathcal N)) \approx \mathbf R$.
\par

For the strong topology, the situation is quite different;
indeed, for ``many'' $\operatorname{II}_1$-factors, for example for those
associated to infinite amenable icc groups or to non-abelian free groups,
it is known that the group
$\operatorname{U}_1(\mathcal N)^{\text{strong topology}}$
is contractible \cite{PoTa--93}.

\section{\textbf{Revisiting the Fuglede-Kadison 
\\
and other determinants}}
\label{FK2}

Most of this section can be found in \cite{HaS--84a}.
For other expositions of part of what follows, 
see \cite[around Theorem 1.10]{CaFM--97}
and \cite[Section 3.2]{Luck--02}.

Let $A$ be a complex Banach algebra (with unit, again for simplicity reasons),
$E$ a Banach space, 
and $\tau : A \longrightarrow E$ a continuous linear map 
which is \emph{tracial},
namely such that $\tau(yx) = \tau(xy)$ for all $x,y \in A$.
Then $\tau$ extends to a continuous linear map 
$\operatorname{M}_\infty(A) \longrightarrow E$,
defined by $x \longmapsto \sum_{j \ge 1} \tau(x_{j,j})$,
and again denoted by $\tau$.
If $e,f \in \operatorname{M}_\infty(A)$ are equivalent idempotents, 
we have $\tau(e) = \tau(f)$;
it follows that $\tau$ induces a homomorphism of abelian groups
\begin{equation*}
\underline{\tau} \, : \,  K_0(A) \longrightarrow E, \hskip.5cm
[e] \longmapsto \tau(e) .
\end{equation*}
For example, if $A = \mathbf C$ 
and $\tau :  \mathbf C \longrightarrow \mathbf C$ is the identity, 
the stable equivalence class of an idempotent $e \in \operatorname{M}_n(\mathbf C)$
is precisely described by the dimension of the image 
$\operatorname{Im}(e) \subset \mathbf C^n$,
so that $K_0(\mathbf C) \approx \mathbf Z$,
and the image of $\underline{\tau}$ is the subgroup $\mathbf Z$
of the additive group $\mathbf C$.
\par

For a piecewise differentiable path $\xi : [\alpha_1, \alpha_2]
\longrightarrow \operatorname{GL}_\infty^0(A)$, we define
\begin{equation}
\label{Eqbasic}
\aligned
\widetilde \Delta_\tau (\xi) \, &= \,
\frac{1}{2 \pi i} \tau \Big( \int_{\alpha_1}^{\alpha_2}
\dot{\xi} (\alpha) \xi(\alpha)^{-1} d\alpha \Big) 
\\
\, &= \, 
\frac{1}{2 \pi i}   \int_{\alpha_1}^{\alpha_2}
\tau \big( \dot{\xi} (\alpha) \xi(\alpha)^{-1} \big) d\alpha .
\endaligned
\end{equation}
(If $X$ is a compact space, 
for example if $X = [\alpha_1, \alpha_2] \subset \mathbf R$, 
the image of a continuous map 
$X \longrightarrow \operatorname{GL}_\infty^0(A)$
is inside $\operatorname{GL}_n(A)$, 
and therefore in the Banach space $\operatorname{M}_n(A)$,
for $n$ large enough; the integral can therefore be defined naively
as a limit of Riemann sums.)
\par

The normalisation in (\ref{Eqbasic}) is such that, 
if $A = \mathbf C$ and $\tau = \operatorname{id}$,
the loop defined by $\xi_0(\alpha) = \exp (2 \pi i \alpha)$ for $\alpha \in [0,1]$
gives rise to $\widetilde \Delta_\tau (\xi_0) = 1$.

\begin{Lemma}
\label{lemmepourtildedelta}
Let $A$ be a complex Banach algebra with unit, $E$ a Banach space,
$\tau : A \longrightarrow E$ a tracial continuous linear map, and
\begin{equation*}
\widetilde \Delta_\tau \, : \,
 \{\text{paths in} \hskip.2cm \operatorname{GL}_\infty^0(A) \hskip.2cm \text{as above} \}
 \longrightarrow E
\end{equation*}
be the mapping defined by (\ref{Eqbasic}).
\begin{itemize}
\item[(i)]
If $\xi$ is the pointwise product of two paths $\xi_1, \xi_2$
from $[\alpha_1, \alpha_2]$ to $\operatorname{GL}_\infty^0(A)$,
then $\widetilde \Delta_\tau (\xi) = \widetilde \Delta_\tau (\xi_1)
+ \widetilde \Delta_\tau (\xi_2)$.
\item[(ii)]
If $\Vert \xi(\alpha) - 1 \Vert < 1$ for all $\alpha \in [\alpha_1, \alpha_2]$,
then $\tau\big( \dot{\xi}(\alpha) \xi (\alpha) ^{-1} \big) d\alpha$
has a primitive $\tau \big( \log \xi(\alpha) \big)$, so that
\begin{equation*}
2 \pi i \widetilde \Delta_\tau (\xi) = \tau( \log \xi(\alpha_2)) - \tau( \log \xi(\alpha_1)) .
\end{equation*}
\item[(iii)]
$\widetilde \Delta_\tau (\xi)$ depends only on the homotopy class of $\xi$.
\item[(iv)] 
Let $e \in \operatorname{M}_\infty(A)$ be an idempotent
and let $\xi_e$ be the loop defined as in (\ref{xie});
then 
\begin{equation*}
\widetilde \Delta_\tau (\xi_e) \, = \,  \tau (e) \in E .
\end{equation*}
\end{itemize}
\end{Lemma}

\noindent
\emph{Sketch of proof.}
Claim (i)  follows from the computation
\begin{equation*}
\aligned
\widetilde \Delta_\tau (\xi_1 \xi_2) \, &= \, 
\frac{1}{2 \pi i}  \int_{\alpha_1}^{\alpha_2}
\tau \Big( \big( \dot{\xi_1} (\alpha) \xi_2 (\alpha) + \xi_1 (\alpha) \dot{\xi_2}(\alpha) \big)
\xi_2(\alpha)^{-1} \xi_1(\alpha)^{-1} \Big) d\alpha
\\
\, &= \, 
\frac{1}{2 \pi i}   \int_{\alpha_1}^{\alpha_2}
\tau \big( \dot{\xi_1} (\alpha) \xi_1(\alpha)^{-1} \big) d\alpha 
\\ &\hskip2cm
\, + \,
\frac{1}{2 \pi i}   \int_{\alpha_1}^{\alpha_2}
\tau \big( \xi_1(\alpha) \dot{\xi_2} (\alpha) \xi_2(\alpha)^{-1} \xi_1(\alpha)^{-1} \big) 
d\alpha 
\\
\, &= \, \widetilde \Delta_\tau (\xi_1) \, + \, \widetilde \Delta_\tau (\xi_2) .
\endaligned
\end{equation*}
Claims (ii) and (iii) are straightforward, 
compare with the end of Subsection \ref{threeformulas}.
Claim (iv) follows again from an easy computation.
\hfill $\square$

\begin{Def}
Let $A$ be a complex Banach algebra with unit, $E$ a Banach space,
and $\tau : A \longrightarrow E$ a tracial continuous linear map.
Define
\begin{equation}
\label{defDelta}
\Delta_\tau : \operatorname{GL}_\infty^0(A) 
\longrightarrow E / \underline{\tau}(K_0(A))
\end{equation}
to be the mapping
which associates to an element $x$ in the domain
the class modulo $\underline{\tau}(K_0(A))$
of $\widetilde \Delta_\tau (\xi)$,
where $\xi$ is any piecewise differentiable path in $\operatorname{GL}_\infty^0(A)$
with origin $1$ and extremity $x$.
\end{Def}

\noindent
\textbf{Remark.}
We insist on the fact that, in general, $\Delta_\tau$ 
is not defined on the whole of $\operatorname{GL}_\infty(A)$.
\par
However, there are several classes of algebras for which is it known
that  the group $\operatorname{GL}_\infty(A)$ is connected.
For example, this is the case for $A = \mathbf C$
(since $\operatorname{GL}_n(\mathbf C)$ is connected for all $n \ge 1$)
and for other finite-dimensional C$^*$-algebras
(which are of the form $\prod_{j=1}^k \operatorname{M}_{n_j}(\mathbf C)$),
more generally for AF C$^*$-algebras,
and also for von Neumann algebras (viewed as Banach algebras).

\begin{Thm} Let the notation be as above.
\begin{itemize}
\item[(i)]
The mapping $\Delta_\tau$ of (\ref{defDelta})
is a homomorphism of groups,
with image $\tau(A) / \underline{\tau}(K_0(A))$;
in particular $\Delta_\tau$ is onto  if $\tau$ is onto.
\item[(ii)]
$\Delta_\tau(e^y)$ is the class of $\tau(y)$ modulo $\underline{\tau}(K_0(A))$
for all $y \in \operatorname{M}_\infty(A)$.
\end{itemize}
\end{Thm}

\begin{Cor}
\label{classicalA}
If $\tau : A \longrightarrow \mathbf C$ is a trace such that
$\underline{\tau}(K_0(A)) = \mathbf Z$,
then 
\begin{equation*}
\exp(2 i \pi \Delta_\tau) : \operatorname{GL}_\infty^0(A) 
\longrightarrow \mathbf C^*
\end{equation*}
is a homomorphism of groups, and
\begin{equation}
\label{detexpA}
\exp(2 i \pi \Delta_\tau) (e^y) \, = \, e^{\tau(y)}
\end{equation}
for all $y \in \operatorname{M}_\infty(A)$ (compare with (\ref{detexpC})).
\par
In particular, 
if $A = \mathbf C$ and if $\tau$ is the identity,
then $\exp(2 i \pi \Delta_\tau)$ is the usual determinant on 
$\operatorname{GL}_\infty(\mathbf C)$.
\end{Cor}

\begin{Cor}
\label{A=M}
If $\mathcal N$ is a factor of type $\operatorname{II}_1$ and $\tau$ its canonical trace,
then $\underline{\tau}(K_0(\mathcal N)) = \mathbf R$,
\begin{equation}
\exp(\operatorname{Re}(2 i \pi \Delta_\tau)) : \operatorname{GL}_\infty(\mathcal N) 
\longrightarrow \mathbf R_+^*
\end{equation}
is a surjective homomorphism of groups,
and its restriction to $\operatorname{GL}_1(\mathcal N)$ is the
Fuglede-Kadison determinant.
\end{Cor}

If $A$ is a \emph{separable} Banach algebra given with\footnote{If 
$A$ is a C$^*$-algebra, any tracial continuous linear form
can be written canonically as a linear combination of four tracial continuous linear forms
which are moreover hermitian ($\tau(x^*) = \overline{\tau(x)}$)
and positive ($\tau(x^* x) \ge 0$ for all $x \in A$).
This is essentially due to Grothendieck 
(\cite{Grot--57}, already quoted in Section \ref{vna});
see \cite[Proposition 2.7]{CuPe--79}.}
a trace $\tau$,
then the range of $\Delta_\tau$ is the quotient of $\mathbf C$
by a \emph{countable} group, by Proposition~\ref{Asep}. 
Suppose more precisely that $A$ is a C$^*$-algebra with unit,
that $\tau$ is a faithful tracial continuous linear form on $A$ 
which is \emph{factorial},
and  that the GNS-representation associated to $\tau$
provides an embedding $A \longrightarrow \mathcal N$
into a factor of type $\operatorname{II}_1$, with $\tau$ on $A$ being the restriction
of the canonical trace on $\mathcal N$.

\begin{Rem}
\label{interp}
Let $A, \tau, \mathcal N$ be as above.
We have a commutative diagram
\begin{equation*}
\begin{array}{ccccccccc}
& &
\operatorname{GL}_\infty(\mathbf C)   & \longrightarrow & 
\operatorname{GL}_\infty^0(A)             &   \longrightarrow &  
\operatorname{GL}_\infty(\mathcal N)  & &
\\
& && & & & & & 
\\
& &
\downarrow 2\pi i \Delta_\tau^{(\mathbf C)} & &
\downarrow  2\pi i \Delta_\tau^{(A)} &  &  
\downarrow  2\pi i \Delta_\tau^{(\mathcal N)} &  & 
\\
& & & & & & & & 
\\
\mathbf C^* & \approx &
\mathbf C / 2 \pi i \mathbf Z & \longrightarrow & 
\mathbf C / 2 \pi i  \underline{\tau}(K_0(A)) & \longrightarrow & 
\mathbf C / 2 \pi i \mathbf R & \approx & \mathbf R_+^*
\end{array}
\end{equation*}
where the top horizontal arrows are inclusions,
the vertical homomorphisms $2 \pi i \Delta$'s are onto,
and the bottom horizontal arrows are onto.
\par
In this sense, $\Delta_\tau^{(A)}$ 
can ve viewed as an interpolation 
between $\Delta_\tau^{(\mathbf C)}$,
which is essentially the classical determinant, 
and $\Delta_\tau^{(\mathcal N)}$,
which is essentially the Fuglede-Kadison determinant.
\end{Rem}

This situation occurs for example if $A = CAR$ (see \ref{afewstandard});
and also if $A = C^*_\lambda(\Gamma)$ is 
the reduced C$^*$-algebra of an icc countable group $\Gamma$.

\subsection{On the sharpness of $\Delta$'s}
\label{subsectionsharp}
Let $A$ be a complex Banach algebra. 
Denote by $E_u$ the Banach space quotient of $A$
by the closed linear span of the commutators $[x,y] = xy-xy$, $x,y \in A$;
thus $E_u = A / \overline{[A,A]}$.
The canonical projection $\tau_u : A \longrightarrow E_u$
is the \emph{universal tracial continuous linear map} on $A$.
In some cases, the space $E_u$ has been characterised:
for a finite von Neumann algebra $\mathcal N$ with centre $\mathcal Z$,
the \emph{universal trace} 
(as defined in \cite[chap.\ III, $\S$~5]{Dixm--57})
induces an isomorphism $E_u \approx \mathcal Z$ \cite[chap.\ 3]{FaHa--80}.
Information on $E_u$ for stable C$^*$-algebras 
and simple AF C$^*$-algebras
can be found in \cite{CuPe--79} and \cite{Fac--82a}.
\par

To the universal $\tau_u$ corresponds the \emph{universal determinant}
\begin{equation*}
\Delta_u : \operatorname{GL}_\infty^0(A) \longrightarrow 
E_u / \underline{\tau_u}(K_0(A)) .
\end{equation*}
Observe that any tracial linear map $\tau : A \longrightarrow \mathbf C$
is the composition $\sigma \tau_u$ of the universal $\tau_u$
with a continuous linear form $\sigma$ on $E_u$.
We have
\begin{equation}
D\operatorname{GL}_\infty^0(A) \, \overset{(1)}{\subset} \,
\ker(\Delta_u) \, \overset{(2)}{\subset}  \,
\bigcap_{\sigma \in (E_u)^*} \ker( \Delta_{\sigma \tau_u} ) \, \subset \, 
\operatorname{GL}_\infty^0 (A) .
\end{equation}
Both $\overset{(1)}{\subset}$ and $\overset{(2)}{\subset}$ can be strict inclusions,
but $ \overset{(2)}{\subset}$ is always an equality if $A$ is separable.
The last but one term on the right need not be closed 
in $\operatorname{GL}_\infty^0(A)$.
For all this, see \cite{HaS--84a}.
\par

Let us agree that the universal determinant is \emph{sharp}
if the inclusions $\overset{(1)}{\subset}$ and $\overset{(2)}{\subset}$ 
are equalities,
equivalently if the natural mapping from the kernel
$\operatorname{GL}_\infty^0(A) / D\operatorname{GL}_\infty^0(A)$
of (\ref{K1etK1top}) to $E_u / \underline{\tau_u}(K_0(A))$
is an isomorphism.
\par

If $A$ is a simple AF C$^*$-algebra with unit, 
its universal determinant is sharp.
More precisely, if $A$ is an AF-algebra with unit,
$\operatorname{GL}_n(A)$ is connected for all $n \ge 1$,
and a fortiori so is $\operatorname{GL}_\infty(A)$.
If $A$ is moreover simple, then 
\begin{equation*}
D\operatorname{GL}_n(A) \, = \,  \ker \left(
\Delta_u : \operatorname{GL}_n(A) \longrightarrow E_u / \underline{\tau_u}(K_0(A))
\right) 
\end{equation*}
for all $n \ge 1$,
and a similar equality holds for $\operatorname{U}_n(A)$ 
and $D\operatorname{U}_n(A)$ \cite[th.\ I and prop.\ 6.7]{HaS--84b}.
If $A$ is a simple C$^*$-algebra with unit which is \emph{infinite},
there are no traces on $A$ \cite{Fac--82a}, and therefore no $\Delta_\tau$,
and $\operatorname{GL}_n^0(A)$ is a perfect group \cite[th.\ III]{HaS--84b}.

Moreover, if $G$ is one of these groups, 
the quotient of $DG$ by its centre is a simple group \cite{HaS--85}.

\section{\textbf{On Whitehead torsion}}
\label{baratintorsion}

We follow \cite{Miln--66}.

\subsection{On $K_1$ and basis of modules} 
Let $\mathcal R$ be a ring; we assume that
free $\mathcal R$-modules of different finite ranks are not isomorphic.
Let $F$ be a free $\mathcal R$-module of finite rank, say $n$;
let $a = (a_1, \hdots, a_n)$ and $b = (b_1, \hdots, b_n)$ be two basis of $F$.
There is a matrix $x \in \operatorname{GL}_n(\mathcal R)$ such that
$a_j = \sum_{k=1}^n x_{j,k} b_k$, 
and therefore a class of $x$ in $\overline{K}_1(\mathcal R)$, denoted by $[b/a]$.
\par

Let 
\begin{equation}
\label{chaincomplex}
C \, : \, 0 \, \longrightarrow \,  
C_n \, \overset{d_n}{\longrightarrow} \, 
C_{n-1} \, \overset{d_{n-1}}{\longrightarrow} \,
\cdots \, \overset{d_2}{\longrightarrow} \, 
C_1 \, \overset{d_1}{\longrightarrow} \, 
C_0 \, \overset{d_0}{\longrightarrow} \, 0
\end{equation}
be a chain complex of free $\mathcal R$-modules of finite ranks
such that the homology groups $H_i$ are also free $\mathcal R$-modules
(the latter is automatic if $H_i = 0$, a case of interest in topology).
Suppose that, for each $i$,  there is given \emph{a basis} $c_i$ of $C_i$, 
and a basis $h_i$ of $H_i$
(the latter is automatic if $H_i = 0$).
\par

At first, assume that each boundary submodule
$B_i$ is also free, with a basis $b_i$.
Using the inclusions $0 \subset B_i \subset Z_i \subset C_i$
and the isomorphisms $Z_i/B_i \approx H_i$, $C_i/Z_i \approx B_{i-1}$, 
there is a natural way to define (up to some choices)
\emph{a second basis} of $C_i$, denoted by $b_ih_ib_{i-1}$.
By definition, the \emph{torsion} of $C$, 
given together with the basis $c_i$ and $h_i$, 
is the element\footnote{The occurence of the same letter $\tau$ for the torsion here
and for traces above has no other reason than standard use.}
\begin{equation}
\label{deftauMilnor}
\tau(C) \, = \, \sum_{i=0}^n (-1)^i
[b_ih_ib_{i-1} / c_i] \, \in \, \overline{K}_1(\mathcal R) .
\end{equation} 
It can be shown to be independent of 
the other choices  made to define $b_ih_ib_{i-1}$;
in particular, the signs $(-1)^i$ are crucial for $\tau(C)$
to be independent of the choices of the basis $b_i$'s.
\par

In case the hypothesis on $B_i$ being free is not fulfilled,
it is easy to check that the $B_i$'s are stably free,
and there is a natural way to extend the definition of $\tau(C)$.
This can be read in \cite[$\S$~1-6]{Miln--66}.
(A $\mathcal R$-module $A$ is \emph{stably free} if there exist
a free $\mathcal R$-module $F$ such that $A \oplus F$ is free.)

\medskip

Suppose now that $C$ is \emph{acyclic}, namely that $H_*(C) = 0$.
There exists a \emph{chain contraction}, namely a degree one morphism
$\delta : C \longrightarrow C$ such that $\delta d + d \delta = 1$,
and therefore an isomorphism
\begin{equation}
\label{d+delta}
d + \delta \vert_{\operatorname{odd}} \, : \,  
C_{\operatorname{odd}} = C_1 \oplus C_3 \oplus \cdots
\, \longrightarrow  \,
C_{\operatorname{even}} = C_0 \oplus C_2 \oplus \cdots .
\end{equation}
Since $C_{\operatorname{odd}}$ and $C_{\operatorname{even}}$ 
have basis (from the $c_i$ 's), this isomorphism 
defines an element in $\overline{K}_1(\mathcal R)$;
we have
\begin{equation}
\label{deftauCohen}
\tau(C) \, = \, \text{class of}  \hskip.2cm 
d + \delta \vert_{\operatorname{odd}} 
\hskip.2cm \text{in} \hskip.2cm
\overline{K}_1(\mathcal R) 
\end{equation}
(see \cite[Chap.~III]{Cohe--73}).
Formula (\ref{deftauCohen}) is sometimes better suited than (\ref{deftauMilnor}).

\subsection{The Whitehead torsion of a pair $(K,L)$}
\label{Wtorsionofpairs}
Consider a pair $(K,L)$ of a finite connected CW-complex $K$
and a subcomplex $L$ which is a deformation retract of $K$;
set $\Gamma = \pi_1(L) \approx \pi_1(K)$.
For a CW-pair $(X,Y)$, consider the complex 
which defines cellular homology theory, with groups
$C^{\operatorname{CW}}_i(X,Y) =
H^{\operatorname{sing}}_i(\vert X^i \cup Y \vert, \vert X^{i-1}\cup Y \vert)$;
here, $H^{\operatorname{sing}}_i$ denotes singular homology with
trivial coefficients $\mathbf Z$,
and $\vert X^i \cup Y \vert$ denotes the space underlying
the union of the $i$th skeleton of $X$ with $Y$.
If $\widetilde K$ and $\widetilde L$ denote the universal covers of $L$ and $K$,
the groups $C^{\operatorname{CW}}_i(\widetilde K,\widetilde L)$ are naturally
free $\mathbf Z [\Gamma]$-modules;
moreover, they have free basis as soon as a choice has been made
of one oriented cell in $\widetilde K$ above each oriented cell in $K$.
For each of these choices, and the corresponding basis,
we have a torsion element
$\tau(C^{\operatorname{CW}}(K,L)^{+\text{choices}}) \in 
\overline{K}_1(\mathbf Z [\Gamma])$.
To obtain an element independent of these choices,
it suffices to consider the quotient
$\operatorname{Wh}(\Gamma) =
K_1(\mathbf Z [\Gamma]) / \langle  \{1,-1\}, \Gamma \rangle$
defined in \ref{K1alg}.
The class 
\begin{equation*}
\tau(K,L) \, \in \, \operatorname{Wh}(\Gamma) 
\end{equation*}
of $\tau(C^{\operatorname{CW}}(K,L)^{+\text{choices}})$ is the
\emph{Whitehead torsion} of the pair $(K,L)$.
In 1966,  it  was known  to be \emph{combinatorially invariant} 
(namely invariant by subdivision of CW-pairs);
more on this in \cite[$\S$~7]{Miln--66}.
Since then, it has been shown to be a \emph{topological invariant}
of the pair $(\vert K \vert, \vert L \vert)$ \cite{Chap--74};
this was a spectacular success of infinite dimensional topology
(manifolds modelled on the Hilbert cube, and all that).

\subsection{On torsion and cobordism}
A \emph{h-cobordism} is a triad $(W; M, M')$
where $W$ is a smooth manifold 
whose boundary is the disjoint union $M \sqcup M'$
of two closed submanifolds,
such that both $M$ and $M'$ are deformation retracts of $W$.
Products $W = M \times [0,1]$ provide trivial examples;
in \cite{Miln--61}, there is a non-trivial example of a h-cobordism
$(W, L \times \mathbf S^4, L' \times \mathbf S^4)$,
with $L$ and $L'$ two $3$-dimensional lens manifolds\footnote{More precisely,
$L$ and $L'$ are quotients of $\mathbf S^3$ 
by free actions of $\mathbf Z /7\mathbf Z$.}
which are homotopically equivalent
but not homeomorphic.
More generally, by a 1965 result of Stallings \cite[$\S$~11]{Miln--66}:
\begin{itemize}
\item[]
\emph{If $\dim M \ge 5$, any $\tau \in \operatorname{Wh}(\pi_1(M))$
is of the form $\tau(W,M)$}
\item[] 
\emph{for some h-cobordism $(W; M, M')$.}
\end{itemize}
Together with the s-cobordism theorem (below), this implies:
\begin{itemize}
\item[]
\emph{For two h-cobordisms $(W_1; M, M_1)$, $(W_2; M, M_2)$} 
\item[]
\emph{such that $\tau(W_1,M) = \tau(W_2,M)$,}
\item[]
\emph{there exist a diffeomorphism $W_1 \longrightarrow W_2$ which preserves $M$.}
\end{itemize}

\par\noindent
A h-cobordism gives rise to a chain complex and a torsion invariant
$\tau(W,M) \in \operatorname{Wh}(\pi_1(M))$,
as in \ref{Wtorsionofpairs}.
Here is the basic \emph{s-cobordism theorem} of Barden, Mazur, and Stallings
\cite{Kerv--65}:
\begin{itemize}
\item[]
\emph{If $\dim W \ge 6$, then $W$ is diffeomorphic to the product $M \times [0,1]$}
\item[]
\emph{if and only if $\tau(W,M) = 0 \in \operatorname{Wh}(\pi_1(M))$.}
\end{itemize}
In particular, 
if $M$ is simply connected, then $W$ 
is always diffeomorphic to $M \times [0,1]$;
this is the \emph{h-cobordism theorem} of \cite{Smal--62}.
\par

For example, if $\Sigma$ is a homotopy sphere 
of dimension $n \ge 6$,
if $W$ is the complement in $\Sigma$ 
of two open discs with disjoint closures, 
and if $S_0, S_1$ are the boundaries of these discs
(they are standard spheres),
then $(W; S_0, S_1)$ is a h-cobordism,
and $W$ is diffeomorphic to $S^{n-1} \times [0,1]$.
It follows that $\Sigma$ is \emph{diffeomorphic}
to a manifold obtained by gluing together the boundaries
of two closed $n$-balls under a suitable diffeomorphism,
and that $\Sigma$ is \emph{homeomorphic}
to the standard $n$-sphere;
the last conclusion is still true in dimension $n=5$.
This is the generalised Poincar\'e conjecture in large dimensions,
established in the early 60's.
The first proof was that of Smale (see \cite{Smal--61}, 
and slightly later \cite{Smal--62}); 
very soon after, 
there have been other proofs of other formulations\footnote{For
a description of the various formulations, written for non-specialists,
see \cite{Miln--11}.} 
of the Poincar\'e conjecture, 
logically independent of Smale's proof but inspired by his work, 
by Stallings (for $n \ge 7$) and Zeeman (for $n \ge 5)$.
The other dimensions were settled much later: by Freedman in 1982 for $n=4$
and by Perelman in 2003 for $n=3$.

\subsection{\textbf{On the Reidemeister-Franz-de Rham torsion}}
\label{RFR}

Since $K_1$ is a functor, any linear representation
$h : \Gamma \longrightarrow \operatorname{GL}_k(\mathbf R)$
provides a ring homomorphism 
$\mathbf Z [\Gamma] \longrightarrow \operatorname{M}_k(\mathbf R)$,
therefore a morphism of abelian group
\begin{equation*}
K_1(\mathbf Z [\Gamma]) \, \longrightarrow \,
K_1(\operatorname{M}_k(\mathbf R)) \, = \, 
K_1(\mathbf R) \, \approx \, \mathbf R^* ,
\end{equation*}
where $\approx$ is induced by the determinant 
$\operatorname{GL}_\infty(\mathbf R) 
\overset{\det}{\longrightarrow} \mathbf R^*$,
and also a morphism $\overline{K}_1(\mathbf Z [\Gamma]) 
\longrightarrow
\overline{K}_1(\mathbf R) \approx \mathbf R_+^*$.
When the representation  is orthogonal,
$h : \Gamma \longrightarrow O(k)$,
this induces a morphism of abelian groups 
$\operatorname{Wh}(\Gamma) \longrightarrow  \mathbf R_+^*$.
\par

For a complex of $\mathbf Z [\Gamma]$-modules $C$ 
with torsion $\tau(C) \in \overline{K}_1(\mathbf Z [\Gamma])$, 
the image of $\tau(C)$
is the \emph{Reidemeister torsion} $\tau_h(C) \in \mathbf R_+^*$,
which is a real number (in fact, $\tau_h(C)$ may be well-defined 
even in cases $\tau(C)$ is not~...).
This is the basic invariant in important work
by Reidemeister, Franz, and de Rham
(earliest papers pulished in 1935).
\par

Given by a Riemannian manifold $M$
and an orthogonal representation $h : \pi_1(M) \longrightarrow O(k)$
of its fundamental group,
one defines a complex $C$ (see (\ref{chaincomplex})) 
of differential forms with values
in a bundle associated to $h$.
Under appropriate hypohesis,
one has a famous \emph{analytical expression}
of the Reidemeister-Franz-de Rham torsion,
and an equality
\begin{equation}
\label{RaySinger}
\aligned
\tau_h(C) \, &= \,  \frac{1}{2} \sum_{k=0}^{n} (-1)^k \ln \det (d_k^*d_k)
\\
\, &= \, \frac{1}{2} \sum_{k=0}^{n} (-1)^k \, k \, \ln \det (d_k^*d_k  + d_{k+1}d_{k+1}^*)
\endaligned
\end{equation}
\cite{RaSi--71}.

\section{\textbf{A few lines on $L^2$-torsion}}
\label{baratinL2torsion}

\subsection{Another extension of ${\det}^{FK}_{\tau}$}
\label{extensionFKL}
Let $\mathcal N$ be a finite von Neumann algebra
and let $\tau : \mathcal N \longrightarrow \mathbf C$
be a finite trace.
For $x \in \mathcal N$, let $(E_\lambda)_{\lambda \ge 0}$
denote the spectral measure of $(x^*x)^{\frac{1}{2}}$.
Define
\begin{equation}
\label{FKL}
{\det}^{FKL}_{\tau}(x) \, = \, \left\{
\aligned
\exp \lim_{\epsilon \to 0+} \int_\epsilon ^\infty \ln (\lambda) \hskip.1cm d \tau(E_\lambda)
\hskip.5cm 
&\text{if} \hskip.2cm \lim_{\epsilon \to 0+} \int_\epsilon ^\infty \cdots > - \infty,
\\
0 \hskip4.5cm
&\text{otherwise.}
\endaligned
\right.
\end{equation}
It is immediate that ${\det}^{FKL}_{\tau}(x) = {\det}^{FK}_{\tau}(x)$
when $x$ is invertible, but the equality does not hold in general
(${\det}^{FK}_{\tau}(x)$ is as in (\ref{analyticextension})). 
For example, if $x \in \operatorname{GL}_1(\mathcal N)$, and 
$X = \begin{pmatrix} x & 0 \\ 0 & 0 \end{pmatrix} \in \operatorname{M}_2(\mathcal N)$,
we have
\begin{equation*}
0 \, = \,  {\det}^{FK}_{\tau}( X ) \, \ne \, 
{\det}^{FKL}_{\tau}( X ) \, = \,  
{\det}^{FKL}_{\tau}( x ) \, = \, 
{\det}^{FK}_{\tau}( x ) \, > \, 0 .
\end{equation*}
The main properties of  ${\det}^{FKL}_{\tau}$,
including
\begin{equation*}
\aligned
& {\det}^{FKL}_{\tau}(xy) \, = \,   {\det}^{FKL}_{\tau}(x)  {\det}^{FKL}_{\tau} (y) 
\hskip.2cm \text{for} \hskip.2cm x,y \in \mathcal N 
\\
& \text{such that} \hskip.2cm 
x \hskip.2cm \text{is injective and} \hskip.2cm
y \hskip.2cm \text{has dense image}
\endaligned
\end{equation*}
are given in \cite[Theorem 3.14]{Luck--02}.

\subsection{On $L^2$-torsion}
Since Atiyah's work on the $L^2$-index theorem \cite{Atiy--76},
we know that (complexes of) $\mathcal N$-modules 
are relevant in topology,
say for $\mathcal N = \mathcal N (\Gamma)$
and for $\Gamma$ the fundamental group of the relevant space.
Let $\mathcal N$ and $\tau$ be as above. Let 
\begin{equation*}
C \, : \, 
0 \, \overset{d_{n+1}}{\longrightarrow} \, C_n \, \overset{d_n}{\longrightarrow}  \, 
C_{n-1} \, \overset{d_{n-1}}{\longrightarrow} \, \cdots  \, \overset{d_2}{\longrightarrow} \, 
C_1 \, \overset{d_1}{\longrightarrow} \, C_0 \, \overset{d_0}{\longrightarrow}  \, 0
\end{equation*}
be a finite complex of $\mathcal N$-modules, 
with appropriate finiteness conditions on the modules 
(they should be projective of finite type),
with a condition of acyclicity on the homology
(the image of $d_j$ should be dense in the kernel of $d_{j-1}$ for all $j$),
and with a non-degeneracy condition on the differentials $d_j$
(which should be of ``determinant class'',
namely $\det_\tau^{FKL}(d^*_{j}d_j)$ should be as in the first case
of (\ref{FKL})).
The \emph{$L^2$-torsion} of $C$ is defined to be 
\begin{equation}
\label{defL2torsion}
\rho^{(2)}(C) \, = \, 
\sum_{k=0}^n  (-1)^k \ln {\det}_{\tau}^{FKL} ((d_j^*d_j)^{1/2})
\, \in \, \{-\infty\} \sqcup \mathbf R  
\end{equation}
(compare with (\ref{RaySinger})).
There is a $L^2$-analogue of (\ref{deftauCohen}),
see \cite[3.3.2]{Luck--02}.

\medskip

$L^2$-torsion, and related notions, have properties which parallel
those of classical torsions, in particular of Whitehead torsion,
and seem to be relevant for geometric problems, e.g.~for understanding 
volumes of hyperbolic manifolds of odd dimensions. 
We refer (once more) to \cite{Luck--02}.

\subsection{Speculation}
It is tempting to ask wether (or even speculate that!) 
modules over reduced C$^*$-algebras $A = C^*_{red}(\Gamma)$
and refinements $\Delta_\tau^{(A)}$
will be relevant one time or another,
rather than modules over $\mathcal N (\Gamma)$
and Fuglede-Kadison determinants ${\det}_\tau^{FK}(\cdot)$. 
Compare with Remark \ref{interp}.

A first and \emph{important} technical problems will be that of extending 
these ``determinants''
$\Delta_\tau^{(A)}$ to singular elements.

\end{document}